\documentclass[9pt]{IEEEtran}
\usepackage{fullpage,psfrag,url,subfigure,algorithm,algpseudocode,amssymb,amsthm,version,color,transparent,epsf,booktabs} 
\usepackage[pdftex]{graphicx}
\usepackage[small,bf]{caption}
\usepackage{pstricks}
\usepackage{authblk}
\usepackage[margin=1.5cm]{geometry}

\usepackage[cmex10]{amsmath}
\newcommand{\BEAS}{\begin{eqnarray*}}
\newcommand{\EEAS}{\end{eqnarray*}}
\newcommand{\BEA}{\begin{eqnarray}}
\newcommand{\EEA}{\end{eqnarray}}
\newcommand{\BEQ}{\begin{equation}}
\newcommand{\EEQ}{\end{equation}}
\newcommand{\BIT}{\begin{itemize}}
\newcommand{\EIT}{\end{itemize}}
\newcommand{\BNUM}{\begin{enumerate}}
\newcommand{\ENUM}{\end{enumerate}}

\newcommand{\BA}{\begin{array}}
\newcommand{\EA}{\end{array}}

\newcommand{\ie}{{\it i.e.}}

\newcommand{\reals}{{\mbox{\bf R}}}

\newcommand{\argmin}{\mathop{\rm argmin}}

{\begin{quote}}{\end{quote}}

\newcounter{oursection}

\newcounter{lecture}

\newcommand\numberthis{\addtocounter{equation}{1}\tag{\theequation}}

\makeatletter
\renewcommand\@makefntext[1]{%
  \noindent\makebox[0em][r]{\@makefnmark}#1}
\makeatother

\title{Dynamic Control and Optimization of Distributed Energy Resources in a Microgrid}

\author[1]{
Trudie Wang \thanks{Mechanical Engineering Department, Stanford University} \
Dan O'Neill \thanks{Electrical Engineering Department, Stanford University} \ 
Haresh Kamath \thanks{Electric Power Research Institute}} 

\date{\today}

\begin{document}
\maketitle

\begin{abstract}

As we transition towards a power grid that is increasingly based on energy from renewable resources like solar and wind, the intelligent control of distributed energy resources (DER) including photovoltaic (PV) arrays, controllable loads, energy storage and plug-in electric vehicles (EVs) will be critical to realizing a power grid that can handle both the variability and unpredictability of renewable energy sources as well as increasing system complexity. In addition to providing added system reliability, DERs acting in coordination can be leveraged to address supply and demand imbalances through demand response (DR) and/or price signals on the electric power grid by enabling continuous bidirectional load balancing. Intelligent control and integration has the capability to reduce or shift demand peaks and improve grid efficiency by displacing the amount of backup generation needed and offsetting the need for spinning reserves and peaking power plants. 

Realizing such a decentralized and dynamic infrastructure will require the ability to solve large scale problems in real-time with hundreds of thousands of DERs simultaneously online. Because of the intractable scale of the optimization problem with variables and constraints for every DER, load and generator online at each time period, we use an iterative decentralized method to operate each DER independently and autonomously within this environment. This method was developed in \cite{KCLB:14} using a distributed algorithm referred to as the Alternating Direction Method of Multipliers (ADMM). Specifically, we consider a commercial site equipped with with on-site PV generation, partially curtailable load, EV charge stations and a stationary battery electric storage (BES) unit for backup. The site operates as a small microgrid that can participate in the wholesale market on the power grid or operate off-grid in an islanded state. 
The ADMM algorithm is deployed within a Model Predictive Control (MPC) framework to allow the microgrid to distribute the optimization among the individual DERs and dynamically adapt to changes in the operating environment while responding to external real-time wholesale prices and potential contingency situations. At each time step, embedded controllers model their own DERs as optimization problems with local objectives subject to individual constraints and forecasts. They then use the ADMM algorithm to solve the problem and obtain a control schedule across the MPC horizon. The local objectives are augmented with a regularization term that includes a simple exchanged message between neighbors in the microgrid. This is the only communication required between DERs. Through the exchange of these messages, the decentralized method rapidly converges to an optimal solution for the entire microgrid when each DER is able to locally solve its own problem efficiently in parallel. Once solved, the controllers execute the first step of the schedule and await the next time step at which point they re-solve the problem using any new information that arrives to augment their forecasts over the planning horizon and account for changes in operating state. This iterative optimization process is repeated for every time step thereafter, ensuring a robust and flexible framework that dynamically adapts to changes in the operating environment. We report results for simulations that demonstrate the ability of this optimization framework to respond dynamically in real-time and minimize total cost to the microgrid while respecting and maintaining the functional requirements of all connected DERs.  

\end{abstract}


\section{Introduction} 



Motivated by environmental concerns, the need to diversify energy sources, energy autonomy and energy efficiency, the penetration of Distributed Generation (DG) from renewable resources like solar and wind is rapidly increasing as the trend moves away from large centralized power stations towards more meshed power transmission on the electricity grid. But as penetration of variable generation sources reaches and exceeds the $10-30\%$ range, matching supply to load will begin to pose a significant challenge using existing centralized dispatch mechanisms \cite{BI:04,THPB:08,DEKM:10,WNRN:08}. With these transformations looming in the near-future, the intelligent integration of Distributed Energy Resources (DER) including photovoltaic (PV) arrays, controllable loads, energy storage and the batteries in plug-in Electric Vehicles (EVs) will become crucial to creating a transaction-based collaborative network that can handle both the variability and unpredictability of renewable energy sources as well as increasing system complexity. While these DERs add system complexity, intelligent control of their power schedules has the potential to serve as a considerable system reliability and stability resource while simultaneously providing a means for greater power system flexibility. This can only be achieved, however, if the ability to solve large scale problems in real-time with hundreds of thousands of devices simultaneously online is integrated into the operation of the power system. Optimization at the distribution level will have to go beyond addressing traditional problems like loss minimization and reactive power compensation and consider the foreseeable transition to a more dynamic power system. DERs can be leveraged to address supply and demand imbalances through Demand Response (DR) and/or price signals on the electric power grid by enabling continuous bidirectional load balancing. DERs with storage capabilities can be particularly useful as they allow the excess output from local generation to be absorbed in situations where the grid cannot. This reduces the need for curtailment and firms up power from DG by providing power during shortfalls. DERs acting in coordination could also help reduce or shift demand peaks and improve grid efficiency by displacing the amount of backup generation needed and offsetting the need for spinning reserves and peaking power plants. The dynamic response of a distributed resource located close to the DG source can effectively act as a buffer to match the availability of generation to the draw from online loads. This is critical in low-voltage networks with a high penetration of DG since the absence of buffering through either DR or storage can result in large voltage variations, uncertainty of power flows, and possibly even reversed power flow which may impact local operation of the grid.

The total installed worldwide capacity for DG systems is expected to grow exponentially over the next $5$ years \cite{EIA:MDG13,L:13,M:13}. Additionally, solar is currently the most rapidly growing generation sector and a significant fraction of this has been due to customer-sited PV generation \cite{EIA:11, IEA:12}. This means that the impact of intelligently managed DERs could be substantial. Rapid control of DERs like storage and controllable loads means that despite limited capacity of the individual unit, the potential aggregate impact of a diverse population of DER units can not only help in operating the grid more efficiently but also provide a significant amount of stability and resiliency at times of sudden increase in demand or loss of generation and during times of fluctuating and intermittent wind or solar power. Particularly advantageous is the fact that many DERs do not require time-critical power. However, use of DERs for grid services also has the potential to delay recovery and present additional strain on an already taxed grid as well as possibly accelerate degradation of the lifetime and functionality of the DER. Such impacts must be balanced against the potential positives of using DERs to improve grid operation. 
Of particular interest is the fleet of electric vehicles starting to come online at commercial facilities as well as in car sharing programs, public transportation and goods delivery services. The primary reason for electrifying the transportation fleet is to substantially reduce the amount of fossil fuels used for transportation. Since most vehicles only travel short distances on any given day, even relatively short range plug-in hybrid EVs can run primarily on electricity. And as the penetration of solar and wind continues to grow rapidly, charging these vehicles will become increasingly clean. With some storage capacity, the EV has also come to represent one of the DERs with the greatest potential for DR. A study conducted by the Electric Power Research Institute (EPRI) using data from the most recent National Household Travel Survey (NHTS) \cite{AD:11} found that across many factors commonly assumed to influence driving behaviors, there is actually little variation in typical driving patterns. This suggests that while EVs will be a significant load on the power grid at the distribution level, they will also be able to provide a consistent resource for energy management services. Plug-in EVs are particularly advantageous since they can be deployed with great flexibility when needed given their relatively long availability window when compared with their actual charge duration. However, as the number of EVs continues to increase, integrating them intelligently will be critical on a power grid that was not inherently designed to handle the additional load. The study concluded that while many of the expected factors like providing high power charging and adding additional charge locations would have little impact on the magnitude of coincident load reduction, load control through optimized charging may be able to significantly increase and improve energy resource utilization.

In this paper, we propose a method to realize the shift from a centrally run power grid to a decentralized network that will enable real-time management and scheduling of EVs at a commercial site with on-site PV generation, partially curtailable load and a battery storage unit for backup. While the primary charging location for the majority of EVs will be at home, workplace charging will be the second most important charging location since the longest vehicle dwell times during the day occur at work \cite{MBK:09,AD:11} . We look at workplace charging because in addition to providing an additional charge location to increase the displacement of gasoline with electricity, it provides useful charging opportunities early in the day when additional load generally has the lowest impact on the grid and it also allows the EVs to directly take advantage of local generation from PV. The regular patterns of arrival and departure at work will also make the vehicles reliable to manage. However, because vehicles arrive at work within a narrower time frame than they arrive at home, there is also potentially a sharp load peak that could occur in the morning at work if the EVs all begin charging immediately. Intelligent load management of the EVs with the help of real-time charge rates or other broadcasted signals could hence be effective if not crucial in flattening load peaks and adding load diversity in addition to allowing the vehicles to take advantage of the PV output. From the utility perspective, this load shifting could have significant benefits for the grid since it would reduce or eliminate the load from vehicle charging coincident with the primary load peak in the evening. The commercial site is connected to the distribution system and operates as a microgrid participating in the wholesale electricity market, with the DERs interacting with each other to conduct arbitrage local generation from the PV array. While each DER is controlled independent of one another, the capacity constraint at the point of common coupling (PCC) means that in addition to minimizing the cost of power to the microgrid, the controllers must cooperate in such a way as to prevent violation of the physical limitations on the line. 
The aim is to show that decentralized, distributed control can increase flexibility and responsiveness to both local and system-wide contingencies. In addition to demonstrating that DER integration leads to load balancing and rapid response times, we also show that DER functionality is maintained and the DERs each operate within their constraints. 

Our method draws upon two key algorithms: using the Alternating Direction Method of Multipliers (ADMM) within a Model Predictive Control (MPC) framework, we can distribute the optimization problem among the independent DERs to dynamically and robustly optimize the power flows between them. While there is an abundance of literature that has looked at the response of DERs at various locations on the grid \cite{YT:12,LCL:11,PSM:10,CN:09}, as well as optimization through MPC for electricity and energy purposes \cite{AMTCT:12, HGL:11,JQBS:11,JJSP:11,NHDH:09,YKTGB:11}, most work has focused on looking at how the DERs can be centrally and remotely controlled. There has only been some recent work done in using distributed dynamic algorithms to optimize the power grid and allow it to operate more efficiently \cite{ALP:11,CCH:11,GTL:11,GTL:12,MCH:13}. These approaches make assumptions about DER behavior that restrict their ability to include new objectives and constraints as well as their adaptability to uncertainty. ADMM allows us to pose the integration of DERs as a completely decentralized control problem whereby each controller only needs to exchanges simple messages with its neighbors in the power network in a relatively unconstrained framework.  The method is iterative, passing messages at every iteration before each DER minimizes its own objective function along with a simple regularization term that depends only on the messages it received in the previous iteration from its neighbors. The only coordination that is required between DERs in the network is synchronizing iterations and the DERs are allowed to model their systems independently and autonomously solve their own local optimization problems with great flexibility. We expand on the ADMM work done previously in our group \cite{BPCPE:11,KCLB:14} by looking at how MPC allows us to extend this framework by integrating prediction and state uncertainties in addition to constraints for real DERs into the distributed algorithms dynamically and robustly. At each time step, the decentralized ADMM algorithm uses this information to jointly minimize the sum of the objective functions of all the DERS in the network over the MPC time horizon subject to local constraints in order to determine the control actions of each DER. The process of acquiring new information, making prediction and optimizing using ADMM is then repeated at the next time step and every time step thereafter. Thus while the decentralized ADMM algorithm distributes computation across all DERs in the microgrid and enables rapid decentralized optimization of dynamic objectives, MPC ensures that the solution is robust to missing information and inaccurate forecasts by reoptimizing at every time step. The microgrid effectively becomes self-correcting and the flexibility in the framework allows DERs to plug into the power network with their own configurations, objectives and constraints specified simply and directly without any reconfiguration since the algorithms are application agnostic and handle various objectives, constraints and forecast models. This is particularly useful at the distribution level where both individual utility and social good will depend on a given operating environment. Each independent DER will be able to act autonomously in a way that enables it to participate responsibly on the grid while still maximizing benefit to itself. In this way, collaboration and cooperation through distributed optimization enables sustainable operation of the grid. DERs can take advantage of dynamic situations where they can work together in order to increase net benefits, ranging from minimizing system costs due to real-time rates to avoiding export of power when there are high peak demand charges or if there is catastrophic failure on the grid. The complexity of managing a decentralized network is divided up among all the controllers to make control of the microgrid tractable and reflective of local needs while leveraging local contribution of distribution-level resources including DG from renewable resources and other controllable DERs. The algorithms can also easily be scaled up with each DER or aggregate of DERs participating at the power system level and helping convert the grid into a flexible and simultaneously more resilient system. Creating a distributed and decentralized system where loads are met locally will additionally result in more efficient use of the transmission and distribution system, allowing for deferment of transmission/distribution line upgrades.

In our model, the power converter at the point of common coupling (PCC) only needs to have access to real-time wholesale prices or know that the grid is in a contingency state and each DER only needs to know its own state and historical power profiles in order to determine control actions for each DER in the microgrid. Using ADMM with MPC to decentralize and distribute the optimization problem, this is sufficient to obtain optimized power profiles that minimize the total objective while maintaining DER functionality and respecting both individual and system constraints. The distributed algorithms also enable us to explore the full potential of EVs connected to the grid. Many vehicle manufacturers and utilities are already considering vehicle-to-grid operation (V$2$G) as a way for vehicles to provide additional grid services that they could potentially be remunerated for \cite{KTv2g:05,DEKM:10,CTD:09}. V$2$G enables bidirectional charging and presents the vehicle batteries with the opportunity to provide various grid services when the vehicles are plugged into the charger. While EVs are unlikely to be able to provide significant base load and peak power services without significant penetration, they have great potential in rapid response short-term power service markets like spinning reserves and regulation \cite{WSE:08,DEKM:10,CTD:09} . 

Strongly related to this is the ability to help smooth the intermittency from renewable sources of energy \cite{KT:05,DEKM:10,CTD:09}. With a bi-directional power converter, plug-in EVs can also provide backup power in situations where the microgrid needs to be isolated from the utility grid, either intentionally due to a fault or other abnormal grid conditions. In this stand-alone mode, the vehicles allow the microgrid to continue operating without power from the grid \cite{CTD:09,PAS:09,KT:05}. However, while there are many efforts looking at how to use EVs to provide additional grid services and make the vehicles more economical and profitable, there is still much uncertainty as to whether the profit provided to each vehicle will be sufficient incentive to the vehicle owner who must give up utility, long-term functionality and privacy of the vehicle in return \cite{DEKM:10,LPMC:10,SE:11}. The algorithms we propose address each of these issues, allowing the vehicle owners to determine their own utility functions, balancing economic benefits against vehicle utility as well as the lifetime degradation from additional cycling on the vehicle. Privacy is also retained since no external control is required to schedule the vehicles. While this means that the system operator cannot know with complete certainty what actions each EV will take, the optimality of the ADMM algorithm does ensure that the aggregate behavior of the microgrid can be aligned with common goals including operational efficiency, flexibility, stability and resiliency. We note here that we deliberately remain agnostic about the viability of V$2$G: the EVs are allowed to both charge and discharge in our model in order to allow the algorithms to determine the tradeoff between competing objectives within each vehicle and from other DERs in the microgrid. The decision as to whether or not it is economical for the EVs to participate in the market for V$2$G will naturally fall out of the solution to the optimization problem. 
In each scenario, we compare avoided costs against potential tradeoffs. We also compare these results to the prescient case for $24$ hour look ahead horizons to demonstrate that even with simple forecasting models, we are able to operate robustly and very near optimal using ADMM with MPC.

The remainder of this paper is organized as follows. In Section \ref{model} we provide the technical details and the formal mathematical definition of our microgrid model, including the objectives and constraints of each DER and the objectives of the microgrid when participating in the wholesale electricity market. Section \ref{method} describes both the ADMM algorithm used to solve the model at each time step as well as the MPC framework that the algorithm operates within. The prediction methods used to make forecasts at each time step are also presented. Section \ref{examples} presents a numerical example using wholesale price schedules taken from the CAISO, PV array and load data taken from a commercial site in Northern California, and EV data from a regional transportation study. The results of our simulations are also given in Section \ref{examples}. Finally, in Section \ref{conclusion} we conclude on our results for distributed dynamic optimization and operation of microgrids with controllable DERs using ADMM and MPC.

\section{Model}
\label{model}

\subsection{System dynamics and constraints}

\begin{figure}
\begin{center}
\includegraphics[width=0.45\textwidth]{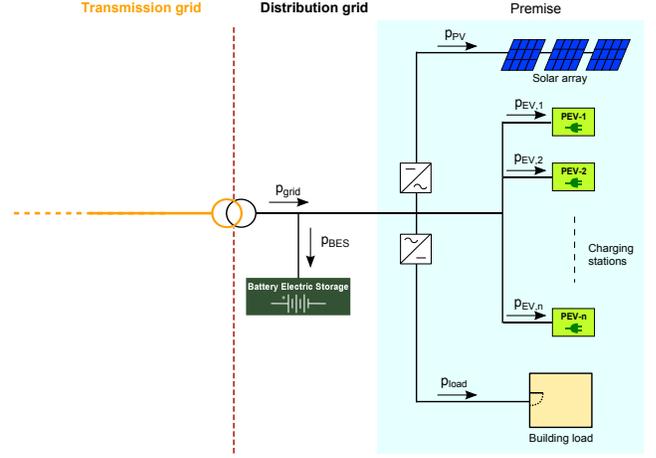}
\end{center}
\caption{Schematic representing microgrid with electric load, PV array, BES and EV charge stations connected to power system. Arrows indicate positive direction of power flow.}
\label{f-model}
\end{figure}

The microgrid modeled and simulated is shown schematically in Fig.~\ref{f-model}. In addition to the electric load onsite and local generation from an installed PV array, there are $N_\mathrm{EV}$ EVs onsite and storage from a BES unit. Collectively, we refer to the PV array, the charging EVs and the BES unit as the DERs. Each DER and the electric load has an associated power schedule across the simulation time horizon $T$, with negative power always defined as being power being generated and/or flowing out from a point in the microgrid. Note that we have also treated the connection at the PCC as an effective DER representing the grid with its own objectives and constraints. 

For simplicity, we consider only DC power without constraints on the phase schedule since smaller microgrids operating at the distribution level will likely be confined to a single phase \cite{bala:08}. The results in this paper can easily be extended to an AC network, however, by imposing an additional phase schedule constraint for each DER modeled \cite{KCLB:14}. We also note here that while all our DERs have been modeled using objective functions that are convex, we do not require either finiteness or strict convexity of any objective function. The convex objective functions and constraints allow the optimization to be carried out rapidly and efficiently, but non-convex problems can also be handled with methods like sequential convex programming and convex relaxation to obtain good local solutions \cite{C:04,KCLB:14,GSC:09,HPR:02}.

The models associated with the microgrid, electric load, and DERs are described in detail below, including the individual objectives and constraints.


\subsubsection{Electric load} 

The electric load is represented with a consumption profile, $p_\mathrm{load}\in\reals^T$. This consumption profile has a diurnal cycle and must be predicted and anticipated using historical load data. We represent the predicted profile as $\hat p_\mathrm{load}\in\reals^T$. When grid-connected or during a contingency situation, the load is curtailable but a minimum time-dependent base load must be met. This is modeled as $\beta \hat p_\mathrm{load} \leq p_\mathrm{load} \leq \hat p_\mathrm{load}$ where $\beta \in [0,1]$ is the minimal load fraction that must be met. The value of $\beta$ will depend on the whether the microgrid is operating in a normally functioning power system or if it is facing a contingency situation.

In order to ensure the amount of load that is curtailed on the system is based on the utility factor of that load, a curtailment cost is also added:
\begin{equation}\label{e-fload}
\begin{array}{l}
\alpha_\mathrm{load} \| \hat p_\mathrm{load} - p_\mathrm{load} \|_2^2
\end{array}
\end{equation}
where $\alpha_\mathrm{load} \geq 0$ is the curtailment penalty parameter.

\subsubsection{Photovoltaic (PV) array}
A PV array also follows a diurnal cycle. The power schedule, $p_\mathrm{PV}\in\reals^T$, is always negative since it only generates power and the peak is offset from the load peak since it occurs when the sun is at its highest. As with the load, the power profile of the PV array is predicted and anticipated using measured historical data. We represent this as $\hat p_\mathrm{PV}$. The array output will fluctuate with cloud cover and contributes to the power balance on the microgrid but the PV inverter has the option to curtail power if it cannot be pushed back out onto the grid or stored. We model this as $p_\mathrm{PV} \leq \hat p_\mathrm{PV}$. 

Since generated power should not be curtailed unnecessarily, a curtailment cost is also included:
\begin{equation}\label{e-fpv}
\begin{array}{l}
\alpha_\mathrm{PV} \| \hat p_\mathrm{PV} - p_\mathrm{PV} \|_2^2
\end{array}
\end{equation}
where $\alpha_\mathrm{PV} \geq 0$ is the PV curtailment penalty parameter.

\subsubsection{Battery Electric Storage (BES).} 
While the inverter of the PV array can do little more than curtail output, a BES unit can store power from the PV array and other DERs and use this to hedge against high prices as well as act as a buffer for requests or unforeseen events on the power grid. Since the BES unit is able to both charge and discharge, its power schedule can be both positive and negative. The rates of charge and discharge are constrained by capacity limits. This is represented as $-D_\mathrm{BES}^\mathrm{max} \leq p_\mathrm{BES} \leq C_\mathrm{BES}^\mathrm{max}$, where $D_\mathrm{BES}^\mathrm{max}$ and $C_\mathrm{BES}^\mathrm{max}$ are the discharging and charging rate limits and $p_\mathrm{BES}\in\reals^T$ is the BES power schedule. 

The dynamics equation governing the state of charge of the BES unit over the time interval $t=1,\ldots,T$ is given by
\begin{equation}\label{e-qbes}
\begin{array}{ll}
q_\mathrm{BES}(t+1) = &\eta_{BES}^q q_\mathrm{BES}(t) + \eta_{BES}^p p_\mathrm{BES}(t),  
\end{array}
\end{equation}
where $\eta_{\mathrm{BES}}^q, \eta_{\mathrm{BES}}^p$ lie in the interval $[0,1]$ and represent the storage and charging efficiencies respectively. The components of $q{BES}\in\reals^{T+1}$ are given as some fraction of the nominal capacity and must remain within the capacity limits of the battery. We specify this as $Q_\mathrm{BES}^\mathrm{min}Q_\mathrm{BES}^\mathrm{cap} \leq q_\mathrm{BES} \leq Q_\mathrm{BES}^\mathrm{max}Q_\mathrm{BES}^\mathrm{cap}$ where $Q_\mathrm{BES}^\mathrm{min}, Q_\mathrm{BES}^\mathrm{max}$ lies in the interval $[0,1]$ and $Q_\mathrm{BES}^\mathrm{cap}$ represents the nominal rating of the unit. The BES unit can also have associated costs, including a cycling cost to penalize excessive charge-discharge cycles, 
\begin{equation}\label{e-fbes}
\begin{array}{l}
\alpha_\mathrm{cyc} \sum_{t=1}^{T-1} |p_{\mathrm{BES}}(t+1) - p_{\mathrm{BES}}(t)|,
\end{array}
\end{equation}
where $\alpha_\mathrm{cyc} \geq 0$ is the cycling penalty parameter.

A terminal constraint can additionally be added to the BES unit to ensure that the storage system is not depleted at the end of the time horizon, with a common choice of $q_\mathrm{final}$ being $0.5Q_{cap}$ or the initial charge state of the battery.

\subsubsection{Electric vehicle (EV)} 

An EV is a flexible load with storage capabilities that allows charging to be deferred as well as discharging when economical. It differs from the BES unit in that it has additional constraints in both availability and required charge capacity. The microgrid has $N_\mathrm{EV}$ EVs onsite and each EV has an associated charging schedule $p_{\mathrm{EV},i} \in \reals^{T}$ and charge state $q_{\mathrm{EV},i} \in \reals^{T+1}$ where $i=1,\ldots,N_\mathrm{EV}$. When the vehicle is not plugged in, $p_{\mathrm{EV},i}(t) = 0$. 
and $t = T_\mathrm{dep}(i)+1, \ldots, T$. The EVs are also associated with four stochastic variables, namely the arrival time, departure time, initial charge state and desired charge state. We write this as a four-component vector $\theta_{\mathrm{EV},i} = [T_\mathrm{arr},T_\mathrm{dep},q_\mathrm{init},q_\mathrm{des}]_{\mathrm{EV},i}$. For each EV, the components of $\theta_{\mathrm{EV},i}$ need to be predicted using probability distributions based on historical data prior to the arrival time. After arrival, the vector is completely defined.

The vehicles also have charge constraints $-C_{\mathrm{EV},i}^\mathrm{max} \leq p_{\mathrm{EV},i} \leq C_{\mathrm{EV},i}^\mathrm{max}$ where $C_{\mathrm{EV},i}^\mathrm{max}$ is the maximum charge/discharge rate of the $i$th EV. 
\begin{equation}\label{e-qev}
\begin{array}{ll}
q_{\mathrm{EV},i}(t+1) = & \eta_{\mathrm{EV},i}^q q_{\mathrm{EV},i}(t) + \\
& \eta_{\mathrm{EV},i}^p p_{\mathrm{EV},i}(t),
\end{array}
\end{equation}
where $\eta_{\mathrm{EV},i}^q, \eta_{\mathrm{EV},i}^p$ lie in the interval $[0,1]$ and are the storage and charging efficiencies respectively.  $q_{\mathrm{EV},i}$ is additionally constrained by the capacity of the battery. This is specified as $Q_{\mathrm{EV},i}^\mathrm{min}Q_{\mathrm{EV},i}^\mathrm{cap} \leq q_{\mathrm{EV},i}(t) \leq Q_{\mathrm{EV},i}^\mathrm{max}Q_{\mathrm{EV},i}^\mathrm{cap}$ where $Q_{\mathrm{EV},i}^\mathrm{min}, Q_{\mathrm{EV},i}^\mathrm{max}$ lie between $[0,1]$ and are the minimum and maximum charge levels of the batteries. $Q_{\mathrm{EV},i}^\mathrm{cap}$ is the nominal rating of the $i$th EV battery.

The desired state of charge for each vehicle $q_{\mathrm{des},i}$ is also used to represent the utility of each vehicle through the constraint
\begin{equation}\label{e-uev}
\begin{array}{l}
q_{\mathrm{EV},i}(T_{\mathrm{dep},i}) \leq \alpha_{\mathrm{des},i} q_{\mathrm{des},i},
\end{array}
\end{equation}
where a higher $q_{\mathrm{des},i}$ indicates a greater demand for power and $\alpha_{\mathrm{des},i}$ represents the vehicle owner's flexibility for not meeting the desired departure state. The utility functions of each vehicle can be learned and/or specified, enabling vehicle owners to assign values to desired charging services and then scheduling their vehicles to maximize net benefits based on these assignments. While methods of determining the utility functions of individual vehicles can readily be incorporated into our model, for simplicity we assume that the utility functions have already been determined and focus on showing that the distributed algorithms can handle varying utility functions for the vehicles simultaneously.

A penalty for excessive cycling of the vehicle battery can also be included for each vehicle,
\begin{equation}\label{e-fev}
\begin{array}{l}
 \alpha_{\mathrm{cyc},i} \sum_{t=1}^{T-1} |p_{{\mathrm{EV},i}}(t+1) - p_{{\mathrm{EV},i}}(t)|,
\end{array}
\end{equation}
where $\alpha_{\mathrm{cyc},i} \geq 0$ is a penalty parameter that weights the excessive cycling cost against the utility function of the vehicle.

\subsubsection{Grid connection} 
The microgrid is connected to the power grid at the PCC where power is stepped up or down in voltage and the wholesale rate for power is applied. From the perspective of the power controller that sits behind the meter, this connection acts as a power limiter that caps the amount of power that can be transmitted over the line. We model this as
\begin{equation}\label{e-ppcc}
|p_{\mathrm{grid}}| \leq P_{\mathrm{PCC}},
\end{equation}
where $p_{\mathrm{grid}} \in \reals^T$ is the power schedule at the PCC and $P_{\mathrm{PCC}}$ is the power limit at the PCC. $p_{\mathrm{grid}}$ can take on both positive and negative values depending on whether the power is being taken from the grid or put back onto the grid. 
At the PCC, the cost function minimizes energy cost as well as the cost of regulation,
\begin{equation} \label{e-fpcc} 
\begin{array}{l}
\sum_{t=1}^T c(t) p_\mathrm{grid}(t)  + 
f_{\mathrm{smooth}}(p_{\mathrm{grid}}(t)).
\end{array}
\end{equation}
In the first term, $c$  is the real-time wholesale price schedule for power. The microgrid receives this price schedule from the independent system operator (ISO) and uses it in the first term to determine what power profile will minimize the cost of power drawn from the grid over the price time horizon. When the term is negative, the microgrid is selling power to the power grid and making a profit. While there are many methods to achieve energy arbitrage and shift demand at the distribution level, the most practical and effective technique is often to simply use real-time wholesale pricing since they are best positioned to allow both consumers and providers to participate and benefit from power transactions as the power grid develops into a more dynamic and responsive system \cite{FS:09,GHSMB:04,HGN:06,BJR:02}. Although such dynamic price schedules are not currently deployed at all levels of the grid, many jurisdictions are moving towards models that pass on prices to consumers reflecting the actual cost of power, especially as the technologies for communicating these prices becomes more common and widespread. 
The second term is the cost associated with the smoothness of the output at the PCC. It ensures a smooth power profile that does not sharply increase or decrease due to PV output or sudden load changes in the microgrid. This measure of power quality can in some cases be remunerated for by the ISO and in other cases is required in order to connect to the power grid and avoid penalties. For this term, we consider a weighted sum of three different measures of smoothness to account for the maximum range, slope and curvature of the power output. This minimizes the variation in magnitude and rate of change of the output to produce a more consistent and smoother power profile. We represent this term as 
\begin{equation} \label{e-fsmooth} 
\begin{array}{l}
f_{\mathrm{smooth}}(p_\mathrm{grid}(t)) = \alpha_{\mathrm{range}} (\max_t p_\mathrm{grid}(t) - \min_t p_\mathrm{grid}(t))\\
\ + \alpha_{\mathrm{diff}} \sum_{t = 1}^{T-1} |p_{\mathrm{grid}}(t+1) - p_{\mathrm{grid}}(t)| \\
\ + \alpha_{\mathrm{curv}} \sum_{t = 1}^{T-2} (p_{\mathrm{grid}}(t) - 2p_{\mathrm{grid}}(t+1) + p_{\mathrm{grid}}(t+2))^2,
\end{array}
\end{equation}
where $\alpha_{\mathrm{range}}$, $\alpha_{\mathrm{diff}}$ and $\alpha_{\mathrm{curv}}$ are the penalty parameters for tuning the range, slope and curvature terms respectively.




\subsection{Objective}

Because our objectives are distinct and separable, we can consider a cost function that is equal to the sum of the individual objectives. This objective is minimized subject to the power flow balance constraint,
\begin{equation}\label{e-fmicrogrid}
\begin{array}{ll}
\mbox{minimize} & \sum_{i = 1}^N f_i(p_i)\\
\mbox{subject to} & \sum_{i = 1}^N p_i = 0,
\end{array}
\end{equation}
%
%
where $N = 4 + N_\mathrm{EV}$ is the total number of independent controllers onsite (\ie load, PV, BES, PCC, EVs) and $f_i$ represents the objective and constraints of each controller. We set $f_i(p_i) = \infty$ to represent infeasibility of the power schedule of the $i$th controller. When $f_i(p_i) < \infty$ for $i = 1,\ldots,N$ and the constraint ensuring power balance is satisfied, the solution is feasible and $p_i$ represents a realizable schedule with $f_i(p_i)$ being the associated cost of that schedule. It is important to note that only the controller at the PCC is provided with the external price schedule and it determines its own power profile. As will be explained in Section \ref{method}, even though it is the interface between the power grid and the other DERs, it does not need to know anything about them or control their behavior. The power balance constraint is a physical constraint that ensures all controllers in the microgrid coordinate to achieve a feasible solution.

\subsection{Control policy}
The control policy for each DER in our microgrid selects the control variables based on information available at the current time. Known information includes DER parameters, grid line limits and conditions at the PCC, and measured states of the DERs and loads. This information along with external wholesale prices and estimates of unknown quantities are then used to calculate and minimize the total cost. Through this optimization process, the control policy determines the power schedules of each controller onsite.

\section{Method} \label{method}

In this section, we describe the algorithms we use to control and optimize the on-site controllers described in Section \ref{model}. Our approach uses the Alternating Direction Method of Multipliers (ADMM) to distribute the optimization process among all the controllers in our microgrid within a Model Predictive Control (MPC) framework. 

Combining ADMM with MPC enables us to take advantage of a symbiotic relationship that enhances the advantages of each. In ADMM, all controllers can compute in parallel and the computation time per iteration is small and independent of the size of the network, enabling MPC to be used in areas where performance was previously limited due to the time required at each computational time step. The computational requirements are effectively reduced by leveraging the local control efforts of each individual controller while respecting the constraints of each DER. And while the complexity of each controller's control effort can remain low, the aggregate intelligence of the coordinated effort is high. The decentralized method of optimization hence addresses the complexity and scale of actively managing a dynamic and rapidly changing power system. But MPC is also a natural extension of ADMM since after one time step is executed, a warm start can be used to rerun ADMM in the subsequent time steps by taking advantage of the power schedules and dual variables computed by each controller in the previous iteration of MPC.  This symbiosis can dramatically speed up computation time to fraction of a second rates \cite{KCLB:14}. This is especially true with recent advances in convex optimization which allow relatively inexpensive embedded processors in power conversion devices to efficiently execute ADMM iterations in tens to hundreds of microsecond time scales \cite{MWB:10, CVXGEN, MWB:11}. 
The combination of a rapid open source convex solver that handles objectives and constraints directly \cite{CVXGEN, RTCVX:10} and increasing CPU capabilities enables power schedules to be dynamically computed. Additionally, offline simulations using different control policies can be benchmarked on similar time scales prior to online implementation to enable design and development of grid systems.  

Another important characteristic that is shared between ADMM and MPC that enables them to work synergistically is that both algorithms only require each controller to have access to \emph{local} variables and not those concerning other controllers in the microgrid.  This allows controllers to operate with their own estimation methods and determine their own control actions without worrying about what other controllers are using since local anomalies are accounted for in the power schedules and messages passed between DERs in ADMM and errors in predictions that may occur at an instance of time are addressed through MPC's iterative process. Also, since all communication is local and dynamic, ADMM used with MPC is ideal for rapidly changing and expanding power networks such as microgrids or commercial sites with energy innovations since it is robust to single points of failure and unexpected topology changes.  Because control is distributed across all DERs, failure of any node or line automatically causes them to adapt and reconfigure their power flows. This ability to self-heal means that the microgrid can handle the increase of DERs in the medium and low voltage distribution networks where the proliferation of high numbers of relatively small individual capacities suggests a new way to think of operating the grid through bottom-up control. Since DG, DR and distributed storage are equivalent in the control sense with an increase in production having the same effect on load balancing as a decrease in consumption, the controller algorithms developed for the microgrid can easily be extended to schedule all types of DERs that connect to a system.

We begin by describing how we can solve the optimization problem described in Section \ref{model} using ADMM equations to distribute the computation effort at a given time step. We then describe the MPC framework that envelopes the ADMM algorithm and allows it to robustly operate and control the DERs in the microgrid, using the ADMM solution to execute local control actions at each time step. Finally, we describe the algorithms we use to make predictions of unknown variables that can then be passed on as inputs to the optimization problem.

\subsection{Alternating Direction Method of Multipliers (ADMM)}
 
In ADMM, we solve the optimization problem specified in \eqref{e-fmicrogrid} using the methods developed in \cite{KCLB:14,BPCPE:11}. The problem is first rewritten in ADMM form by making copies of the variables: 
\begin{equation}\label{e-admm}
\begin{array}{ll}
	\mbox{minimize} & \sum_{i = 1}^N f_i(p_i) + g(\sum_{i = 1}^N z_i)\\
	\mbox{subject to} & p_i - z_i = 0 \quad i = 1,\ldots,N,
\end{array}
\end{equation}
with variables $p_i,~z_i \in \reals^T$. $f_i$ forms the local cost function for controller $i$, and the shared objective $g$ is a function of the sum of the power variables. By introducing $z_i$~'s into the problem as copies of the power variables, the shared constraint represented by the power flow balance can be moved into $g$ which simply becomes the indicator function of the empty set $\emptyset$, \ie is represented by the sum of the variables. 
This definition of $g$ leads to a special case of the sharing problem referred to as the optimal exchange problem. In exchange ADMM, we split the additive term into separable objectives which can then be updated separately to drive the variables to consensus. Each controller in the microgrid effectively participates in an internal market with a price adjustment process that is used to attain general market equilibrium. The internal price of power is increased or decreased depending on whether there is an excess demand or excess supply respectively. This internal price will naturally reflect the external price at the PCC but is not necessarily equivalent since it also reflects the individual interests of the independent DER controllers within the microgrid. The complete derivation of ADMM for the exchange problem can be found in \cite{KCLB:14,BPCPE:11}. We highlight the main equations from the derivation below and ask the reader to refer to \cite{KCLB:14,BPCPE:11} for more details. 

ADMM can be derived directly from the augmented Lagrangian
\begin{equation*}
\begin{array}{ll}
L_\rho(p,z,y) &= \sum_{i = 1}^N \left( f_i(p_i) + g\left(\sum_{i = 1}^N z_i\right) \right. \\
& \left. + y_i^T (p_i-z_i) + (\rho/2) \|p_i-z_i\|_2^2 \right),
\end{array}
\end{equation*}
where $\rho > 0$ is the \emph{penalty parameter}. Minimizing the variables and updating the dual variable independently and iteratively gives the ADMM algorithm
%
%
\begin{eqnarray}
p_i^{k+1} &:=& \argmin_{p_i} \left(f_i(p_i) + (\rho/2)\|p_i - z_i^k + u_i^k\|_2^2 \right) \label{e-admm-p-update-alt} \\
z^{k+1} &:=&  \argmin_{z} \left(g(\sum_{i = 1}^N z_i) \right.\\
 & & \left. + (\rho/2)\sum_{i = 1}^N\|z_i - p_i^{k+1} - u_i^k\|_2^2 \right) \label{e-admm-z-update-alt}\\
u_i^{k+1} &:=&  u_i^k + p_i^{k+1} - z_i^{k+1}. \label{e-admm-u-update}
\end{eqnarray}
with the penalty parameter $\rho$ in the Lagrangian used as the step size in the dual variable update in order to ensure the iterate $(p_i^{k+1},y_i^{k+1})$ is dual feasible \cite{BPCPE:11}. We have combined the linear and quadratic terms in the augmented Lagrangian in order to present ADMM in a more concise and convenient form using a scaled form of the dual variable. The variables $p$ and $z$ are minimized in step~\eqref{e-admm-p-update-alt} and step~\eqref{e-admm-z-update-alt} respectively, followed by a dual variable update in~\eqref{e-admm-u-update}. Written in this form, we can also clearly interpret the scaled dual variable as being the running sum of the residuals,
%
%
%
We can further simplify the $z$-update step by taking the Lagrangian and solving the dual problem to determine $z$ analytically \cite{BPCPE:11}. Substituting the solution in the $u$-update shows that the dual variables $u_i^k$ are all equal and can be replaced with a single dual variable $u$ representing the scaled price. 
The ADMM algorithm finally becomes
%
\begin{eqnarray}
p_i^{k+1} &:=&  \argmin_{p_i} \left( f_i(p_i) + (\rho/2)\|p_i-p_i^k+\overline p^k+u^k\|^2_2\right) \label{e-admm-p-update-fin}\\
u^{k+1} &:=&  u^k + \overline p^{k+1}. \label{e-admm-u-update-fin}
\end{eqnarray}
where $\overline{p} = (1/N) \sum_{i = 1}^N p_i$ represents the mean of the $p_i$ variables. 

%

In each iteration of the $p$-update in step \eqref{e-admm-p-update-fin}, ADMM augments its own local objective function $f_i$ with a simple quadratic regularization term. The linear parts of the quadratic terms containing the iterative target value are then updated, pulling the variables towards an optimal value and allowing them to converge.  
This proximal operator contains the scaled price $u^k$ and can be interpreted as a penalty for $p_i^{k+1}$ deviating from $p_i^k$ projected onto the feasible set, helping pull the variables $p_i$ toward schedules that enable power balancing while still attempting to minimize each local objective. In other words, it represents each controller's commitment to help reach market equilibrium so that as the power profiles are adjusted and the system converges, the effect of the proximal regularization term vanishes. Since each controller only handles its own objectives and constraints, the $p$-update can be carried out independently in parallel by all the controllers in the microgrid. ADMM enables decentralization by allowing the proximal regularization term in the augmented Lagrangian to separate so that it can be minimized locally while still maintaining the convergence properties and robustness of the method of multipliers \cite{BPCPE:11}. This means that the optimization problem with at least as many variables as the number of DERs and loads multiplied by the length of the time horizon is reduced to small local optimization problems with only a few variables for each controller to manage, including local power flows and internal states. The controllers pass their updated power schedules, $p_i^{k+1}$, to a collector (possibly co-located at the PCC) in the $u$-update \eqref{e-admm-u-update-fin} which in turn simply gathers the variables and computes the new average power imbalance, $\overline{p}^{k+1}$, in order to update the scaled price $u^{k+1}$. The computed values are then broadcasted back to the controllers to readjust the proximal operator. In this way, the collection stage projects the power schedules back to feasibility and helps push the system towards equilibrium by adjusting the price up or down depending on whether there is net power demand or generation in the system. The iterative ADMM algorithm converges by alternating between the controllers and the collector with synchronization (necessary for real-time pricing in any event) being the only coordination that is required between the controllers. 

When implementing ADMM, selecting the correct penalty parameter $\rho$ is critically important for convergence rates. The optimal value of $\rho$ will greatly depend on the scheduling problem. While there are heuristic methods to help determine the value of $\rho$, in many cases it will perform just as well with a fixed value found using a binary search method. Since we are using a scaled form of ADMM, the scaled dual variable $u^k = (1/\rho)y^k$ must also be rescaled after updating $\rho$. This means that if $\rho$ is halved, $u^k$ should be doubled before computing the ADMM updates.

It is instructive to consider the primal and dual residuals in the ADMM exchange problem since both have useful interpretations. Following from \cite{BPCPE:11,KCLB:14}, the primal and dual residuals for the exchange problem can simply be written as
\[
r^k= \overline p^{k+1}, \qquad s^k = \rho((p^k-\overline{p}^k)-(p^{k-1}-\overline{p}^{k-1})).
\]
The primal residual simply represents the net power imbalance across all controllers and is hence a measure of physical feasibility. The dual residual is equal to the difference between the current and previous iterations of the deviation in power schedules. As the residuals approach $0$, the internal price represented by the dual variable $y$ converges to an optimal value. 
While there can be more than one optimal point, ADMM is guaranteed to converge to an optimal point when the problem is closed, convex and proper \cite{BPCPE:11}. The residuals also provide a simple criterion for terminating the ADMM algorithm,
\[
\begin{array}{c}
\|r^k\|_2 = \|\overline{p}^k\|_2 \leq \epsilon^\mathrm{pri}, \\ 
\|s^k\|_2 = \|\rho((p^k-\overline{p}^k)-(p^{k-1}-\overline{p}^{k-1}))\|_2 \leq \epsilon^\mathrm{dual},
\end{array}
\]
where $\epsilon^\mathrm{pri}$ and $\epsilon^\mathrm{dual}$ are the primal and dual tolerances respectively. 

\paragraph{Discussion.}
We briefly summarize the ADMM exchange algorithm here in the context of the power flow model. ADMM solves the power flow problem by distributing computational requirements across the multiple controllers. This casts the optimization problem as a completely decentralized control problem whereby each controller computes and exchanges simple proximal messages with only its neighbors in the microgrid. The controllers send small quantities of numeric data to neighbors in order to coordinate at each iteration while storing small amounts of state information and efficiently computing solutions for its own local optimization. In this way, the optimization problem is solved locally in a peer-to-peer fashion and the computational requirements to solve the problem are significantly reduced. While we do not explore the rates of convergence and scalability using ADMM in this paper, we direct the reader to \cite{KCLB:14} for an extensive look at the algorithm's ability to operate in real-time on large scale systems with minimal computational requirements.

External price signals can still be used as an input to help the power grid achieve system-wide objectives but each local controller exchanges the proximal messages which can be thought of as an \emph{internal} price signal that all controllers agree on in order to align locally optimized operating policies with the goals that benefit the entire system. At convergence, this optimal \emph{internal} price between participating DERs in the microgrid naturally results from running ADMM and represents the equilibrium price that occurs when the objectives are mutually optimized. 
Since optimization is independent and allows for autonomous operation with minimal coordination, this bottom-up control approach not only reduces the communication requirements but also makes it feasible to connect large numbers of these distributed systems to the grid without requiring the implementation of complex top-down control systems that require extensive empirical knowledge of each DER. This paradigm shift allows for a new way to think of operating the grid since ADMM can allow for both efficient energy trade and active flexible control of power flows at the controller level. The cooperative approach between all the controllers searches for an outcome that each finds at least acceptable and for which the total objective function representing the social benefit to all participants is minimized.

\subsection{Model Predictive Control (MPC)}

Model predictive control (MPC) is a control policy that can be used to dynamically control each controllable DER independently onsite. Its ability to handle uncertainty and dynamics in real-time allows ADMM to operate robustly in a changing environment. Within the MPC framework, the optimization problem is solved dynamically at each time step using the decentralized ADMM algorithm to determine an action policy for each controller in the microgrid and enable scheduling of the controllable DERs over a finite time horizon. In order to solve the optimization problem, predictions of local load and PV output based on historical data are required. Although a formal statistical or stochastic model can be used to represent uncertainty when making predictions in MPC, it is not needed for the policy to work and the controller can often perform very robustly even when predictions are poor \cite{WB:08}.
  
To implement MPC, we first solve the optimization problem at the current time step $t$ using ADMM to determine a series of conditional power schedules for each controllable DER over a fixed time horizon extending $T_\mathrm{MPC}$ steps into the future. The controllers use the information along with predictions of unknown quantities available to them locally at that time to jointly minimize their objective functions over the MPC time horizon subject to individual constraints. After a prediction is made, ADMM messages are passed between controllers to enable the distributed algorithm to converge rapidly to a solution. Since convergence is on the order of milliseconds and microseconds, actions and schedules are determined well before the next time step in MPC at which the microgrid responds to the updated state of the power system (generally on the order of minutes). Each controller then executes the first step of its schedule and idles until the next time step at which point the entire optimization process is repeated in order to incorporate changes in operating environment, as well as new state measurements and external information that may have subsequently become available.  The process of acquiring new information, making predictions and optimizing using ADMM is repeated at every time step thereafter. This iterative process using MPC effectively ensures that the ADMM solution is robust to measurement errors, missing information and inaccurate forecasts, allowing the controllers to adjust their schedules in response to external disturbances that were unknown at the original time the schedules were computed and ensures that the control policy dynamically adjusts and is self-correcting as new information arrives and changes in the operating environment occur. MPC is thus well suited for use with ADMM in dynamic operation of a microgrid, especially when there is uncertainty in the system \cite{JMM:02, Bem:06, MWB:11}. 

Algorithm~\ref{alg-mpc} outlines the iterative process used by the MPC method. 

\begin{algorithm}
  \begin{algorithmic}
			\State Initialize charge states $q_\mathrm{BES}(1), q_{\mathrm{EV},i}(1) \quad \forall i = 1,\ldots,n$ 
			\For{$t=1 \to T$}
   			\State \parbox[t]{\dimexpr\linewidth-\algorithmicindent}{Predict $\hat p_\mathrm{PV}$ and $\hat p_\mathrm{load}$ using updated historical power profiles available at time $t$;\strut}
				\State \parbox[t]{\dimexpr\linewidth-\algorithmicindent}{Update EV parameters $\hat \theta_\mathrm{EV}$ using charge data as well as arrival and departure times available at time $t$;\strut} 
				\vspace{0.02in}
				\State Solve ADMM problem over MPC time horizon $[t, t+T_\mathrm{MPC}]$; 
				\vspace{0.02in}
				\State $q_\mathrm{BES}(t+1) \leftarrow \hat q_\mathrm{BES}(t+1)$;
				\State $q_\mathrm{EV,i}(t+1) \leftarrow \hat q_\mathrm{EV,i}(t+1) \quad \forall i = 1,\ldots,n$; 
			\EndFor
	\end{algorithmic}
\caption{Iterative optimization using MPC.}
\label{alg-mpc}
\end{algorithm}

MPC is self-adjusting and since forecasting is carried out at every time step with parameters refitted using updated information, a model that is able to represent the general dynamics of the DERs and load is adequate. MPC also directly integrates objectives and constraints without having to learn and adjust controller parameters via a trial and error process, and this makes it particularly well suited for rapid real-time optimization of a locally controlled system which is dependent on and particular to the operating environment. The fact that both our objective and constraints are convex means the problem can be solved very efficiently and while MPC is a heuristic policy that is generally not optimal, it often performs far superior to traditional control methods \cite{MRRS:00,BV:09,CVX}. 

In order to avoid oscillations between MPC iterations, we have also added a regularization term to each objective function $f_i$ which we specify as $f_i^\mathrm{prev} = \alpha^\mathrm{prev} ||\hat p_i - \hat p_i^\mathrm{prev}||_2^2$, where $\hat p_i^\mathrm{prev}$ is the solution for the previous iteration and $\alpha^\mathrm{prev}$ is the damping weight for oscillations between iterations. This prevents the solution at the current time step from diverging too far from the previous time step and ensures that smoothness across the MPC prediction horizon does not come at a cost to smoothness between the current time step and the previous time steps.


%

\section{Predictions} \label{predictions}

To implement MPC with ADMM, estimates of input variables are required at each time step in order to solve the decentralized optimization problem and determine a control policy over the finite time horizon. These estimates can be based on historical data, stochastic models, forecasts, and pricing information. The flexibility of the MPC framework and the decentralized nature of ADMM means that they are not tied to any specific forecasting method and can incorporate a range of techniques depending on what information can be accessed. Perfect predictions or a formal statistical model to represent uncertainty are also not required for the method to perform robustly and predictions that capture general trends are sufficient since MPC recalculates power schedules at each time step after executing the first step of the schedule determined by the ADMM algorithm, dynamically adjusting and self-correcting for any past errors or missing information (\cite{WB:08}).

\subsection{Power profiles}

Within the microgrid, predictions are required for the output of the PV array and the electric load. While many methods have been used to predict energy profiles on the power grid, ranging from multiple regression to expert systems \cite{AN:02, DKO:12, BK:01} and total load over a large region can be predicted up to $1\%$ accuracy, predicting generation locally from intermittent renewable resources poses a greater challenge due to nuances that require detailed knowledge of the system environment and geographic diversification cannot smooth out \cite{FERC:12}. Since our objective is to only capture general trends and the DERs just needs to forecast their own future power profiles in each time period, we require only local measurements of historical output to get adequate predictions. 

To predict PV output, we use the prediction model defined in \cite{WKW:14}. We first assume that the historical PV output data $p_\mathrm{PV}^\mathrm{hist}\in\reals^T_\mathrm{hist}$ has general periodicity over a $24$hour period, where $T_\mathrm{hist}$ is the historical time horizon. This is a good assumption over a period of a few days where the seasonal variation is insignificant since the solar insolation at any geographical location and given time is well determined. Deviations from the expected PV array output are due to weather conditions like cloud cover that result mostly in drops in output with occasional over-irradiance due to cloud enhancement effects. This leads up to propose an asymmetric least squares fit of the available historical PV output data to provide a periodic baseline that is weighted towards the outer envelope of the observed data. To generate the PV output prediction, we first determine the historical baseline $\hat p_\mathrm{PV}^\mathrm{hist}$ by solving an approximation problem with a smoothing regularization term and a periodicity constraint,

%
\begin{align*} 
\underset{\hat p_{\mathrm{PV}}^\mathrm{hist}}{\mbox{minimize}} \ & \frac{1}{T}\sum_{\tau=t-T_\mathrm{hist}+1}^{t-1} \left((\hat p_{\mathrm{PV}}^\mathrm{hist}(\tau) - p_\mathrm{PV}^\mathrm{hist}(\tau))_+^2 \numberthis \label{e-predict-baseline} \right.\\ 
& + \gamma_{\mathrm{asym}} (\hat p_{\mathrm{PV}}^\mathrm{hist}(\tau) - p_{\mathrm{PV}}^\mathrm{hist}(\tau))_-^2 \\ 
& \left.\vphantom{((\hat p_{\mathrm{PV}}^\mathrm{hist}(\tau) - p_\mathrm{PV}^\mathrm{hist}(\tau))_+^2} + \gamma_{\mathrm{curv}} (\hat p_{\mathrm{PV}}^\mathrm{hist}(\tau-1) - 2\hat p_{\mathrm{PV}}^\mathrm{hist}(\tau) + \hat p_{\mathrm{PV}}^\mathrm{hist}(\tau+1))^2\right)    \\
\mbox{subject to} \ & \hat p_{\mathrm{PV}}^\mathrm{hist}(\tau) = \hat p_{\mathrm{PV}}^\mathrm{hist}(\tau + T_\mathrm{period}) \\
& \qquad \qquad \tau=t-T_\mathrm{hist}+1,\ldots,t-1,  
\end{align*}
where $(z)_+ = \max(0,z)$ and $(z)_- = \min(0,z)$. $\gamma_{\mathrm{asym}}$ and $\gamma_{\mathrm{curv}}$ represent the weights for the asymmetric and curvature terms respectively and $T_\mathrm{hist}$ is the amount of time over which $p_{\mathrm{PV}}^\mathrm{hist}$ occurs. 
The first and second term of the objective function represent the positive and negative deviation of the predicted curve from the actual data and the third term smooths the curvature in the predicted curve. The first term is weighted more heavily to push the predicted curve towards the outer envelope of the fluctuating data. The constraint ensures periodicity across $T_\mathrm{period} = 24\mathrm{hrs}$. Solving this problem de-noises the data to reconstruct a smooth baseline profile \cite{BV:09}. The objective function is a weighted sum of squared convex terms and forms a regularized convex problem which trades off an asymmetric least-squares fit against the mean-square curvature of the data. The baseline prediction for the MPC horizon is then defined as $\hat p_{\mathrm{PV}} = \hat p_{\mathrm{PV}}^\mathrm{hist}(t-T_\mathrm{period},\ldots,t+T-T_\mathrm{period})$. 

Once we determine this baseline prediction, we correct for transient weather phenomena by adjusting the baseline using an error fit with a linear model applied to the residual $r = \hat p_{\mathrm{PV}}^\mathrm{hist} - p_{\mathrm{PV}}^\mathrm{hist}$. The correction is calculated by writing the predicted residual at the time step $\tau$ as a weighted sum of the previous residuals, $\hat r(\tau) = a_1r_{\tau-1} + a_2r_{\tau-2} + \cdots + a_nr_{\tau-n}$ where $a$ is an $n$-element vector determining what weights to give the $n$ previous residuals. The associated residual prediction error is defined as $e(\tau) = \hat r(\tau) - r(\tau)$. To determine $a$, we minimize the sum of the squared $l_2$ norms of the prediction error over the entire historical data time horizon $\tau=t-T_\mathrm{hist}+1,\ldots,t-1$. 
We can rewrite this concisely in matrix form, 
\begin{equation}\label{e-errfit}
||e||_2^2 = ||Ma - b||_2^2
\end{equation}
where
\resizebox{0.9\columnwidth}{!}{
$
b = \left[\BA{c} r(n+1) \\ r(n+2) \\ \vdots \\ r(T_\mathrm{hist}) \EA \right], \quad
M = \left[ \BA{cccc} r(n) & r(n-1) & \vdots & r(1) \\ r(n+1) & r(n) & \vdots & r(2) \\ \vdots & \vdots & \vdots & \vdots \\
r(T_\mathrm{hist}-1) & r(T_\mathrm{hist}-2) & \vdots & r(T_\mathrm{hist}-n)  \EA \right] 
$
}\vspace{0.1in}

This is a least squares problem and has the analytical solution $a = M^\dagger b$ where the ${}^\dagger$ symbol denotes the Moore-Penrose pseudoinverse. We assume $M$ is full rank since the residuals are due to random weather patterns and can effectively be treated as independent and identically distributed.
The predicted residual corrections across the MPC horizon are then decreased by a factor $\lambda$ at each future time step, where $0 < \lambda < 1$. This reduces the magnitude of the correction over the MPC horizon moving forward in time so that the prediction reverts back to the baseline. The procedure used to make predictions at each MPC time step is summarized in Algorithm~\ref{alg-prediction}. 

\begin{algorithm}
  \begin{algorithmic}[1]
		\State Compute baseline profile $\hat p_\mathrm{PV}^\mathrm{hist}$ over MPC horizon using $p_\mathrm{PV}^\mathrm{hist}$ to solve \eqref{e-predict-baseline}.
		\State Compute residual $r = \hat p_{\mathrm{PV}}^\mathrm{hist} - p_{\mathrm{PV}}^\mathrm{hist}$.
		\State Determine the residual weights $a$ by \ref{e-errfit}.
		\State Compute predicted residuals over MPC horizon as $\hat r(\tau) = \lambda^{\tau-t} a^T[r(\tau-1) \cdots r(\tau-n)]$ for $\tau = t,\ldots,t+T$.
		\State Compute prediction as $\hat p_{\mathrm{PV}} + \hat r$.
\end{algorithmic}
\caption{Computing PV predictions from historical data.}
\label{alg-prediction}
\end{algorithm}


We employ a similar approach for the electric load using the historical load profile $p_\mathrm{load}^\mathrm{hist}\in\reals^T_\mathrm{hist}$. The only difference between the PV and load profiles is that with the load, we weight the first and second term of the objective function equally since we expect the positive and negative deviation of the predicted curve from the actual data to be similar. The adaptability of MPC means that the least squares fit of the load output and the adjusted least squares fit of the PV output data is able to provide sufficiently accurate forecasts to enable dynamic power scheduling. This prediction method is simple to implement using available historical data and requires very little computational effort during real-time implementation. 

\begin{figure}
	\begin{center}
	\subfigure{
		\centering
		\includegraphics[width=0.45\textwidth]{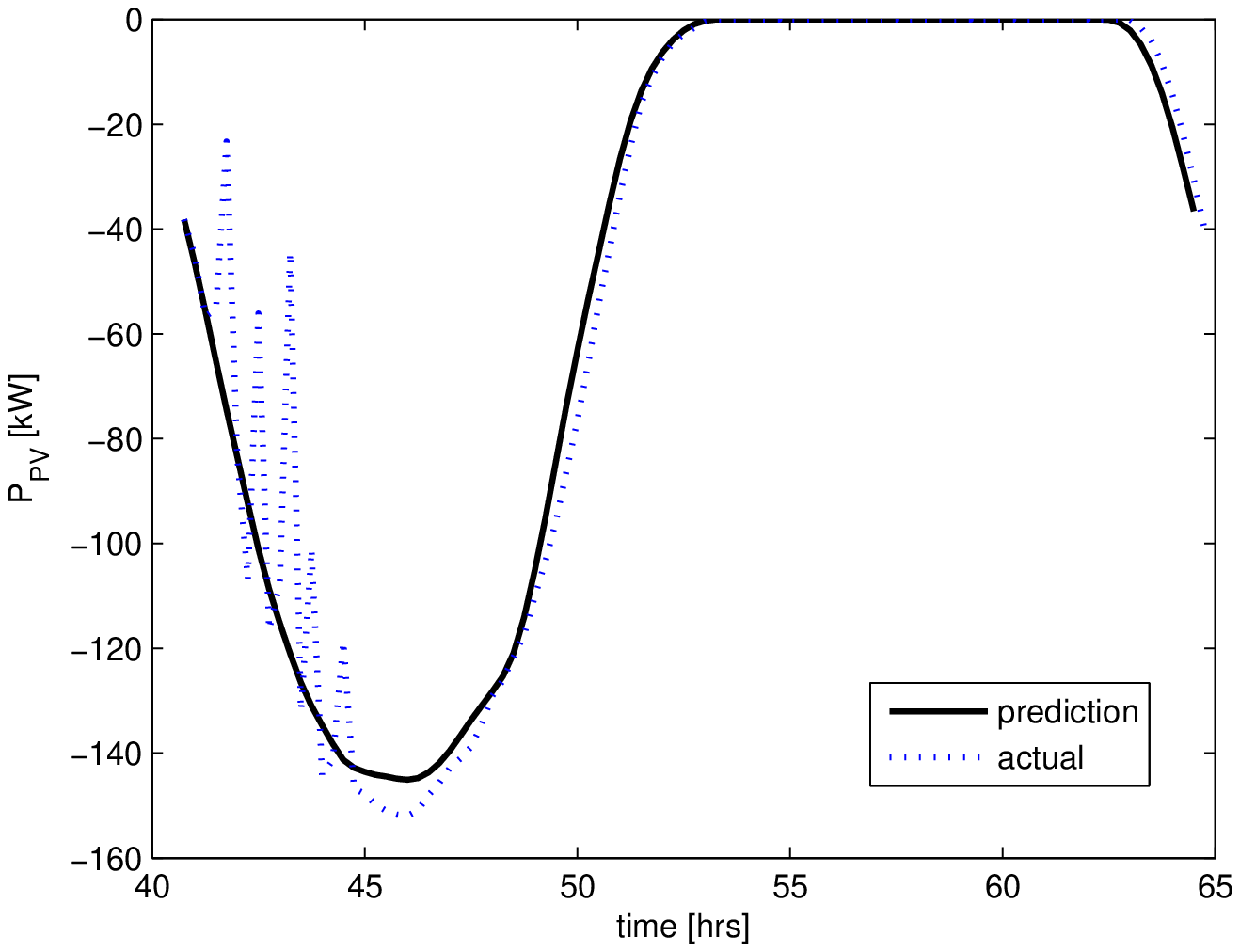}
	}
	\hspace{0.75cm}
	\subfigure{
		\centering
		\includegraphics[width=0.45\textwidth]{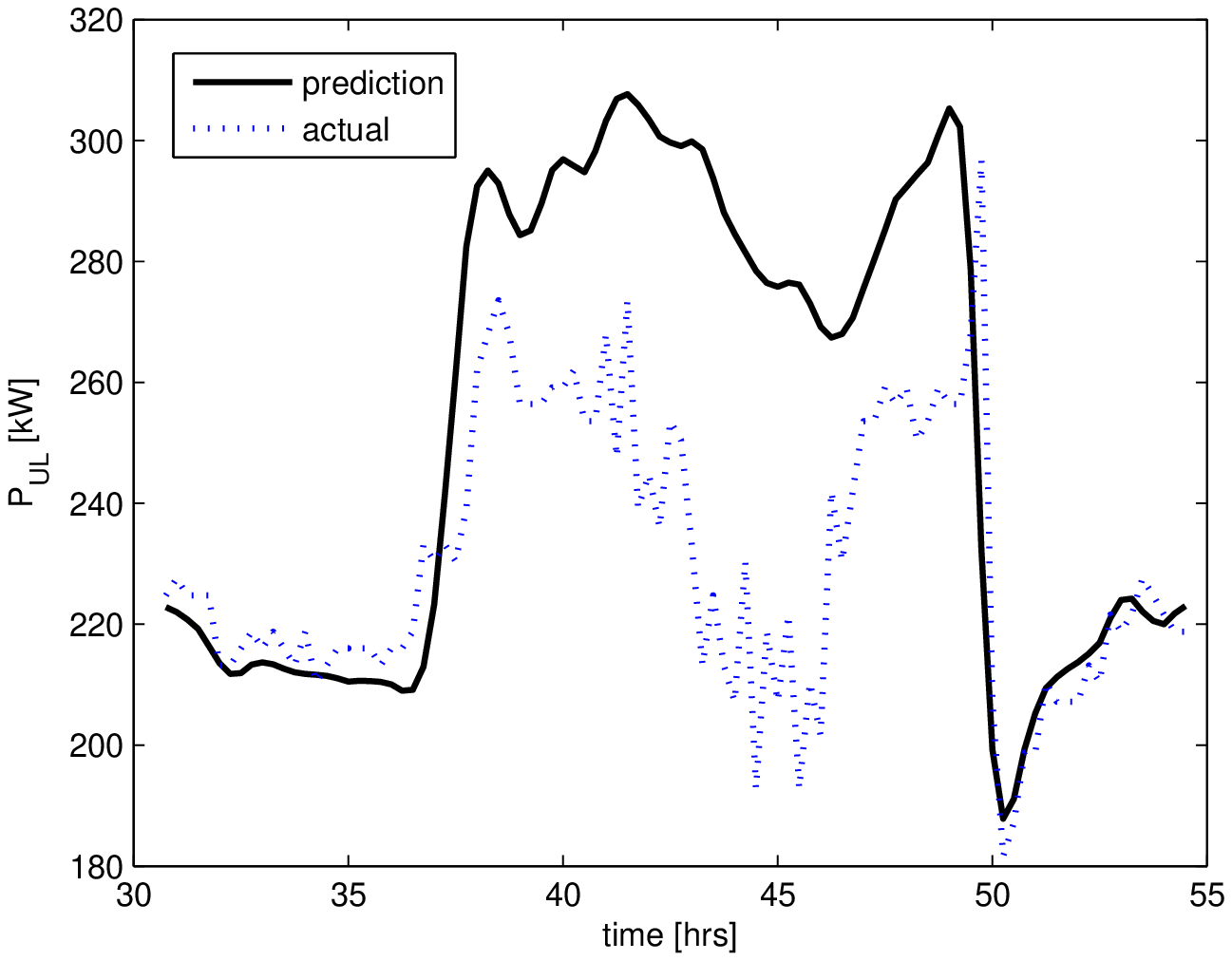}
	}
	\end{center}
\caption{Predictions for the PV output (\emph{top}) and electric load profile (\emph{bottom}) over the MPC time horizon at one instance of time.}
\label{f-predict}
\end{figure}

Examples of PV and load profile predictions over the MPC time horizon at one instance of time are presented in Fig.~\ref{f-predict}. In each plot, the prediction is compared to the actual output. For both PV and load, the predictions are good at capturing the diurnal trends and general shape of the power profiles. Sudden changes in power profiles that occur near the time of measurement where the actual output is known and rapid deviations can be anticipated are also captured by the error correction, pulling the baseline towards the actual output and ensuring immediate response to sudden changes. Additional measurements and information to supplement the historical power profiles used to make the predictions could help to capture the irregular and intermittent dips and spikes that occur throughout the prediction horizon. This additional information might include cloud cover predictions for PV or anticipating occupancy levels in the case of load. While the flexibility of our prediction method means that it is trivial to incorporate this added data by including an additional weighted term to the asymmetric least squares objective, MPC does not require more than what is shown in Fig.~\ref{f-predict} in order to perform robustly. 

\subsection{EV parameters}
In order to predict the stochastic variables in $\theta_{\mathrm{EV},i}$ for $i=1,\ldots,n$, we first consider whether or not a given EV has arrived. In the case where an EV has arrived, we assume $\theta_{\mathrm{EV},i}$ is known. If a vehicle has not arrived, we need to predict $\theta_{\mathrm{EV},i}$ using available data.  Ideally, there would be enough data to form approximately stationary prior probability distributions for each of the variables in $\theta_{\mathrm{EV},i}$. In practice, there will initially be insufficient data to make accurate predictions for each vehicle. As time progresses, however, accumulation of data will lead to more stationary distributions and predictions will become increasingly accurate. In fact, for this reason it is common practice to ignore some number of samples at the beginning before the stationary distribution to be reached. However, for MPC even a general forecast can produce good results as our results will later show. 



To predict $T_\mathrm{arr}$ and $T_\mathrm{dep}$,  we use maximum likelihood estimates from stationary probability distributions. While we could construct conditional probability distributions that are time dependent using Markov chains that adjust the stationary prior probability distribution, we found that this method added insignificant or no benefit since MPC already corrects for prediction error and hence already performs near optimal with rough predictions. The fact that the stationary probability distributions are fairly time independent after enough data has been acquired means that this prediction method performs satisfactorily and is preferred for its minimal computational requirements.

Since there is no reason to assume the state of charge variables are strongly and independently correlated with the time variables, we similarly use stationary distributions to determine arrival and desired departure charge states. For  $q_\mathrm{init}$, this means computing the maximum likelihood once daily for the vehicles that have not yet arrived. To determine the desired charge state $q_\mathrm{des}$, we employ a more conservative approach to determine how much charge the vehicle batteries have to at least have at departure time. Instead of a maximum likelihood, we consider the highest charge state in the distribution since it is assumed that the driver will want to ensure there is enough charge to face the majority of contingency situations. 
Although this is overly conservative in the case of a plugin hybrid EV which has the option to use gasoline as backup, gasoline prices will inevitably continue to trend upward and running solely on the electric battery will become increasingly desirable so that the assumptions made serve as a good first approximation for setting a desired charge state. From a carbon standpoint, ensuring the maximum number of electric miles will also have the best possible societal benefit. As more data for each vehicle is accumulated, the model can easily be adapted to the needs of individual vehicle owners in real-time. After a vehicle arrives, the charge states along with the times of arrival and departure are added to the data and the predictions are updated the next day.


\section{Numerical example} \label{examples}

Simulations are run using real load and generation data taken from May~$18-25,2013$. For the simulation period, we also use day-ahead hourly wholesale price data published by the CAISO \cite{CAISO}. In each scenario we simulated, we selected an MPC time horizon of $T=96$. This corresponds to $15$ minute intervals over a $24$ hour period and is a typical horizon and time step for schedule updates. The time step $\tau=1$ corresponds to midnight.


The microgrid is connected to the distribution system through a bidirectional meter which has access to the real-time wholesale prices. The connection also has a physical transfer constraint $P_{\mathrm{PCC}} = 200\mathrm{kW}$. The power profile for the electric load is taken from measurements at a typical commercial site and the PV output is taken from the output of a rooftop PV array, both geographically co-located in Northern California. We size the array to $1200$~kW and scale the output accordingly so that it is able to meet all of the local energy demands of the load over a diurnal cycle when islanded. This data is incorporated into the model both as simulation data and as historical data to help make predictions and schedule the DERs. As time proceeds, the simulation data is added onto the historical data incrementally at each time step and used to update the predictions. The predictions are made using $5$ days of historical data to predict the next $24$ hours of generation. While it may seem that $P_{\mathrm{PCC}}$ is overly limiting, the aim is to demonstrate that the microgrid can run flexibly and reliably on a grid of limited capacity using dynamic algorithms to handle the coupling constraints. The on-site BES unit is sized to provide sufficient capacity to arbitrage energy for islanding situations while maintaining at least $50\%$ baseload \cite{CDCLH:11,CGW:12}. All BES parameters are provided in table \ref{t-BES}. 

\begin{table}
	\centering
	\begin{tabular}{ l | c }
		Parameter & Value \\
		\hline\hline
		$C_\mathrm{BES}^\mathrm{max} = D_\mathrm{BES}^\mathrm{max}$ & $500\mathrm{kW}$\\
		$\eta_{BES,c} = \eta_{BES,d}$ & $0.85$\\
		$\eta_{BES,q}$ & $0.90$\\
		$Q_\mathrm{BES}^\mathrm{min}$ & $0.20$\\
		$Q_\mathrm{BES}^\mathrm{max}$ & $0.90$\\
		$Q_\mathrm{BES}^\mathrm{cap}$ & $3000\mathrm{kWh}$
	\end{tabular}
	\caption{BES parameters.}
	\label{t-BES}
\end{table}


Within the microgrid, a fleet of $20$ EVs is available and capable of level~$2$ charging at $7.2$kW without requiring a dedicated circuit \cite{DOE-EV:13}. Battery efficiency values $\eta_{\mathrm{EV},i}^q, \eta_{\mathrm{EV},i}^p$ are both set at $90\%$ \cite{CHD:10}. The battery capacity for each vehicle is selected based on driver need. In other words, as an appropriate first approximation it is assumed that vehicles with longer commute distances will have vehicle owners who desire larger batteries. Using individual vehicle data taken over several hundred days in a study conducted by EPRI \cite{alex:12}, we consider the longest trip distance and use a $2\times$ buffer along with a typical mileage conversion rate of $0.311$kWh/mile \cite{AD:11}. The minimum and maximum charge states of the batteries are selected to be $30\%$ and $90\%$ respectively to avoid deep cycling of the battery \cite{CHD:10}. The arrival times $T_\mathrm{arr}$ and departure times $T_\mathrm{dep}$ for each vehicle are chosen from distributions constructed using the vehicle data. The value of $q_\mathrm{init}$ is determined using a Monte Carlo method to select from the distribution of initial charge states. For the desired charge state $q_\mathrm{des}$, we use the conservative approach described in Section~\ref{predictions}. 

\section{Results} \label{results}

\begin{center}
\begin{table*}[ht]
{\small
\hfill{}
 \begin{tabular}{ l | l | l | l | l | c }
     Method & External cost [\$] & Smoothing cost [kW] & PV curtailed [kWh] & Load curtailed [kWh] & EV shortfall [kWh] \\
     \hline\hline
     centralized (prescient) & $-1682.59$ & $14.99$ & $902.63$ & $122.05$ & $0$ \\
     centralized & $-1680.62$ & $16.36$ & $918.48$ & $114.94$ & $0$ \\
     ADMM decentralized & $-1674.74$ & $16.67$ & $919.64$ & $116.65$ & $0$
 \end{tabular}
}
\hfill{}
\caption{Cost metrics for simulated scenarios.}
\label{t-results}
\end{table*}
\end{center}



The results of our simulations over a $3$ day period 
are provided below in Table~\ref{t-results}. In addition to carrying out the simulation using the ADMM method to distribute the optimization calculations among the controllers, we simulate a case with a microgrid carrying out the same calculations through centralized optimization as well as a centralized case with prescient knowledge in order to provide a benchmark for how ADMM with MPC performs even with simple prediction methods. In the prescient case, instead of using predictions we assume PV and load schedules as well as EV parameters are fully known over the MPC horizon. For each scenario simulated, we calculate the cost of energy at the PCC as well as for each DER. We also determine the ramping cost at the PCC represented by \eqref{e-fsmooth}, the power curtailed by the PV array and load, and the total energy shortfall when charging the EVs to their desired departure state of charge. 

\begin{figure}
	\begin{center}
		\includegraphics[width=0.3\textwidth]{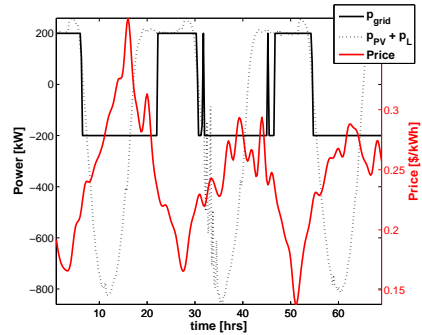}
	\end{center}
\caption{Power profile of microgrid at PCC over the simulation time horizon using ADMM. The output is plotted against price and the unoptimized net power profile of the PV array and load.} 
\label{f-grid}
\end{figure}

\begin{figure}
	\begin{center}
		\includegraphics[width=0.3\textwidth]{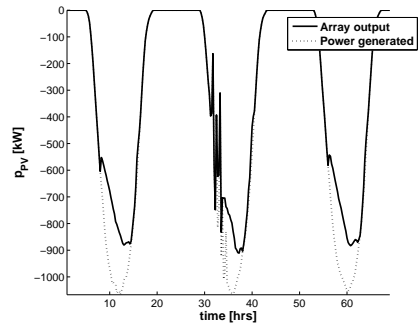}
	\end{center}
\caption{Power output of the PV array over the simulation time horizon using ADMM. The output is plotted against the actual power produced by the array before curtailing.}
\label{f-pv}
\end{figure}

\begin{figure}
	\begin{center}
		\includegraphics[width=0.3\textwidth]{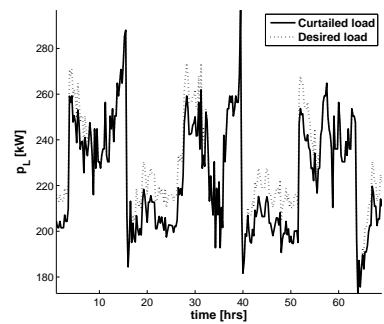}
	\end{center}
\caption{Power profile of the curtailed building load over the simulation time horizon using ADMM. The output is plotted against the desired power prior to curtailing.}
\label{f-ul}
\end{figure}

\begin{figure}
	\begin{center}
		\includegraphics[width=0.32\textwidth]{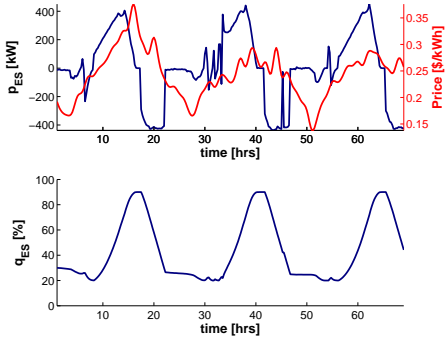}
	\end{center}
\caption{The \emph{top} plot shows the power profile of the BES unit plotted against price using ADMM. In the \emph{bottom} plot, the state of charge of the BES unit is plotted as a percentage of the total battery capacity.}
\label{f-es}
\end{figure}

\begin{figure}
	\begin{center}
		\includegraphics[width=0.32\textwidth]{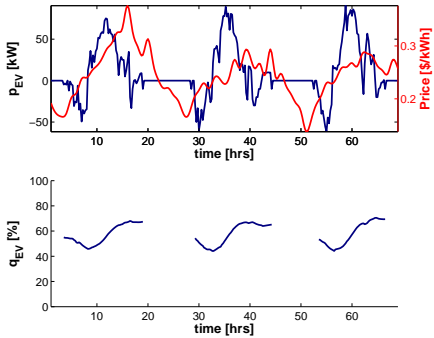}
	\end{center}
\caption{The \emph{top} plot shows the net power profile of the EVs plotted against price using ADMM. In the \emph{bottom} plot, the state of charge of the EVs is plotted as a percentage of the total vehicle battery capacity.}
\label{f-ev}
\end{figure}

As the table shows, ADMM performs similarly to the centralized method and both come very close to the performance of the centralized prescient method. Comparing the ADMM and the centralized case, there is less than a $1.5\%$ difference between the total system costs. Contrasting both of them against the ideal centralized prescient case, there is less than a $3.3\%$ difference between the total system costs. Below, plots of the ADMM and the centralized cases for the simulation period are provided. 
The net power profile of the microgrid at the PCC is shown first in Fig.~\ref{f-grid},  plotted against the net power profile that results from PV and load without any optimization to compare the performance of the microgrid against the base case. Note that we have moved $p_{\mathrm{grid}}$ to the other side of the power balance constraint in the scheduling problem in order to maintain the convention that consumed power is positive and generated power is negative. This plot is followed by plots of the PV array output in Fig.~\ref{f-pv} and the load power profile in Fig.~\ref{f-ul} after curtailing. The power profile and state of charge of the BES unit and the EVs are provided in Fig.~\ref{f-es} and Fig.~\ref{f-ev} respectively. The state of charge in the latter case is plotted as a percentage of the total capacity of all the vehicles beginning each day at the time when the first vehicle arrives and ending at when the last vehicle departs. From the plots comparing ADMM with the centralized case, it can be seen that the results are very similar and this corroborates the quantitative results from Table~\ref{t-results}.

Since the capacity limit of the line is binding, the DERs use their resources to ensure the microgrid can still operate within these limits, absorbing excess generation from the PV array to sustain the load when PV generation drops off. This can be seen in Fig.~\ref{f-grid}, where some smoothing in the power profile is noted with fewer intermittent spikes due to either load or generation but the PCC profile is shaped primarily by the capacity limit. The generated output is relatively flat throughout the peak price hours before dropping off each evening as the price begins to decrease. There is then a constant load for the remaining hours until the PV output begins to increase again along with the price. The PV output is initially insufficient to meet the load which begins to go up earlier so both the BES unit and EVs discharge slightly to help with the load in the face of increasing prices. There is still some flexibility with the storage unit and EVs after the PV output falls off each day, allowing the DERs to use their capacity for energy arbitrage and respond to the still high prices by shifting the power generated by the PV array from its peak time at $12$PM to $5$PM when the peak price occurs. This completely offsets the load at the peak price time in order to increase net profit.  

From Fig.~\ref{f-pv}, there is evidence of some curtailing of generated PV power during peak generation hours when there is not enough capacity both on-site and through the grid connection to absorb that energy. The power is curtailed around the peak generation hour when the price for power is still relatively low. In Fig.~\ref{f-ul}, a small amount of load curtailment also occurs during the evening hours when no solar energy is being generated and the line capacity has been reached in terms of how much power can be imported. In both cases, however, the curtailing is not substantial since both the BES unit and the EVs use their storage capabilities to minimize these effects. The BES unit undergoes a deep discharge daily as shown in Fig.~\ref{f-es} during peak prices when there is no PV output and charges again to full capacity when prices drop and PV output increases. The EVs behave similarly in Fig.~\ref{f-ev}, with the aggregate of vehicles charging as the price generally decreases and discharge excess capacity when the price generally increases. The reason for the smooth charging profile of the vehicles is that while individual charging can be more intermittent when responding to real-time pricing and each vehicle has arrival and departure constraints that limit their flexibility, the aggregate effect of a diverse population can actually produce a smoother net power profile with some degree of coordination without having to impose a smoothing constraint or centralized top-down control \cite{OPMO:13,WGM:13}. This is helped by the fact that when optimizing power schedules through ADMM, the line limit forces each EV to account for the implications of simultaneously drawing power at start up when voltages are low and this presents a natural way to reintroduce diversity by automatically giving priority to the vehicles most willing to accept the high costs without any subjective ordering. Hence the vehicles can be viewed within the microgrid as one single entity with predictable emergent characteristics that can be incentivized to charge in a way which helps to meet system wide objectives even though each vehicle is controlled independently in a decentralized fashion and does not need to over-cycle significantly on an individual basis. This is a notable finding since there is currently a large amount of concern over how the proliferation of EVs may tax an already strained grid with their relatively substantial charging requirements at the distribution level \cite{bullis:13}. Both the predictability of the net power profile and its flexibility to charge at specific hours shows in fact that the EVs can play a very useful role in DR even when control is distributed to each vehicle. Responding in this way through the demand side can provide a potentially more economical way to automate and dynamically respond to changing conditions on the power grid compared to transformers, load tap changers and capacitor banks. Note also that while only the controller at the PCC is aware of the wholesale price schedule and the grid connection limit, all DERs on-site operating independently can help shape $P_\mathrm{grid}$ through ADMM and ensure the limiting constraints are not violated.

While the primary aim of the simulations was to demonstrate the ability of dynamic distributed algorithms to schedule DERs cooperating when faced with a coupling constraint, this example also shows that a microgrid using such algorithms can produce very reliable and consistent power profiles when provided with the right incentives even if its only generation source is intermittent solar. This means that power can be provided with higher power quality and a higher load factor in addition to responding to real-time prices. Such services are critical in a more distributed and decentralized grid where contingency events can lead to catastrophic circumstances as the grid becomes increasingly less foreseeable and controllable. The ability to shape $p_\mathrm{grid}$ using the capacities and demands of each DER benefits both the system operator and the microgrid since it means the system can be run more efficiently at a lower cost and the impacts of tiered pricing can also be avoided.

\section{Conclusion} \label{conclusion}
We have developed a framework using MPC and ADMM to control and optimize a microgrid with PV, curtailable load, EV charge stations and a stationary BES unit. Our work extends previous work done in \cite{KCLB:14} that focused on using ADMM to solve the power scheduling problem. When used with MPC, the algorithms distribute and decentralize the optimization problem and require only local information and simple prediction methods to work while retaining the ability to incorporate any additionally available information into the objectives and constraints of both the entire system and individual DERs as time progresses. By distributing control through ADMM, the problem of managing the microgrid is made more tractable by enabling a cooperative approach between the resources while still respecting coupling constraints due to capacity limited lines. MPC ensures the executed power schedules adapt and keep the system both flexible and resilient when there is imperfect information or sudden unexpected changes to the system. Using data taken from the the EPRI campus in Northern California and transportation data taken from a national survey, we simulated and compared the performance of the algorithms. Our simulations demonstrate that each device is able to retain functionality while allowing the microgrid to respond to external price signals and physical power line limits as well as contingency events. With minimal information sharing between devices, we can obtain performance results that are comparable to those obtained when the optimization problem is solved centrally with prescient knowledge.

\section*{Acknowledgments}
The authors thank Eric Chu and Matt Kraning for extensive discussions on the problem formulation and prediction methods and Stephen Boyd for his advising and guidance on the problem formulation and prediction methods. 

\newpage
\bibliography{smart_converters}

\end{document}